\documentclass{article}

\usepackage{amsmath}
\usepackage[latin1]{inputenc}
\usepackage{graphicx}
\usepackage[T1]{fontenc}
\usepackage{pstricks}
\usepackage{amsthm}

\newtheorem{defi}{Definition}[section]
\newtheorem{theo}[defi]{Theorem}
\newtheorem{cor}[defi]{Corollary}
 
\newtheorem{lemme}[defi]{Lemma}

\newcommand{\R}{\mathbf {R}}
\newcommand{\N}{\mathbf N}
\newcommand{\D}{\mathbf D}
\newcommand{\hd}{\hat{d}}
\newcommand{\acc}{^}
\newcommand{\e}{\varepsilon}
\newcommand{\en}{\varepsilon_n}

\title{A discrete version and stability of Brunn Minkowski inequality }
\author{Michel Bonnefont \\ 
{\small Institut de math\'ematiques, Laboratoire de Statistique et 
Probabilit\'es,} \\ {\small Universit\'e Paul Sabatier,} \\ {\small 118 
route de Narbonne,} \\ {\small 31062 Toulouse,} \\ {\small FRANCE}}
\date{\today}
\begin{document}

\maketitle

\begin{abstract}
In the first part of the paper, we define an approximated Brunn-Minkowski inequality  which generalizes the classical one for length spaces. Our new definition based only on distance properties allows us also to deal with discrete spaces. 
Then we show   the stability of our new inequality under a convergence of metric measure spaces. This result gives as a corollary  the stability of the classical Brunn-Minkowski inequality for geodesic spaces. 
The proof of this stability was done for different inequalities (curvature dimension inequality, metric contraction property) but as far as we know not for the Brunn-Minkowski one.\\
In the second part of the paper, we show that every metric measure space satisfying classical Brunn-Minkowski inequality can be approximated by discrete spaces with some approximated Brunn-Minkowski inequalities. 
\end{abstract}

\section{Introduction}
Let us recall some facts about the Brunn-Minkowski inequality. First the inequality was set in $\R^n$ for convex bodies by Brunn and Minkowski in 1887 (for more details about the inequality and its birth, one can refer to the great surveys \cite{Barthe,Gardner} and the reference therein). It can be  read as if $K$ and $L$ are convex bodies (compact convex sets with non empty interior)  of $\R^n$  and $0<t<1$ then 
\begin{equation}\label{bm1}
V_n((1-t)K+ t L)^{1/n}\geq (1-t) V_n(K)^{1/n} + t V_n(L)^{1/n}
\end{equation} 
 where $V_n$ is the Lebesgue measure on $\R^n$ and $+$ the Minkowski sum  which is given by
 $$A+B= \{ a+b, a\in A, b\in B \}
 $$
 for $A$ and $B$ two sets of $\R^n$. Equality holds if  and only if $K$ and $L$ are equals up to translation and dilatation.
 
 Brunn-Minkowski inequality  is  a very powerful inequality with a lot of applications. For example it implies very quickly the isoperimetric inequality for convex bodies in $\R^n$ which reads
\begin{equation}\label{isop}
\left(\frac{V_n(K) } {V_n(B)} \right) ^{1/n} \leq \left(\frac{s(K) } {s(B)} \right) ^{1/(n-1)}
\end{equation} 
 where $K$ is a convex body of $\R^n$ and $s$ the surfacic measure, with equality if and only if $K$ is a ball.
 
 The Brunn-Minkowski inequality is not only true for convex bodies but also for all compact sets  and  even  for all measurable sets of $\R^n$ (with the little difficulty that the Minkowski sum of two mesurable sets is not necessary measurable). One way to prove it is to prove a functional inequality known as Prekopa-Leindler inequality which applied to characteristic functions of sets gives the multiplicative Brunn-Minkowski
 inequality 
 \begin{equation}\label{bminf}
 V_n((1-t)K+ t L)\geq V_n(K)^{1-t}   V_n(L)^{t}
 \end{equation}
 where $V_n$ is the Lebesgue measure on $\R^n$, $K$ and $L$ two measurable sets of $\R^n$. By homogenity of the volume $V_n$, it can be shown that this a priori weak inequality is in fact equivalent to the $n$-dimensional one (\ref{bm1}).

All this was to show that Brunn-Minkowski inequality has a very geometric meaning and it is natural to ask on which more general spaces than $\R^n$ the inequality can be extended.

One first answer is we can change the measure, for example a measure log-concave on $\R^n$ satisfy multiplicative Brunn Minkowski.

But  to be able to quit $\R^n$, we have to generalize the Minkowski sum. This can be done on length spaces by using ideas of optimal transportation (refer to \cite {Burago} for length space, \cite {Villani} for optimal transporation, and for exemple \cite{CorderoMCS} for this generalisation). Following an idea of this paper,  for two sets $K$ and $L$ of a metric space $X$ 
we define what we are going to call the $s$-intermediate set between  $K$ and $L$ by
\begin{equation}
Z_s(K,L)= \left\{ z\in X; \exists (k,l)\in K\times L, 
\begin{array}{ccc}        
d(k,z)&=&s d(k,l)\\
d(z,l)&=& (1-s) d(k,l)
\end{array}
\right\}
\end{equation}
This set will play the role set of  barycenters of the Minkowski sum.  
 In fact the authors in \cite{CorderoMCS} use it only for a Riemannian manifold but it makes sense for all metric spaces even if it  is interesting only for length space. In this context we will say a metric measure space $(X,d,m)$ satisfies  the $N$-dimensionnal Brunn-Minkowski inequality  if 
 \begin{equation} \label{bmls}
m^{1/N}({Z_s(K,L)})\geq (1-s)\, m^{1/N}(K)  + s \,  m^{1/N}(L)
\end{equation} 
 for all $0<s<1$ and $K$, $L$ compacts of $X$. We will refer in the sequel at (\ref{bmls}) as the "classical"  $N$-dimensionnal Brunn-Minkowski inequality.
 It is proven in \cite{CorderoMCS} that  a  Riemannian manifold $M$ of dimension $n$ whose Ricci's curvature is always non negative satisfies (\ref{bmls}) with dimension $N=n$ and with the canonical volume of the Riemannian manifold as measure, i.e.
  \begin{equation}\label{bm2}
vol(Z_s(K,L))^{1/n}\geq (1-s) vol(K)^{1/n} + s\, vol(L)^{1/n}
\end{equation}  
for all compacts $K$, $L$ of $M$ where $vol$ denotes the canonical volume of the Riemannian manifold. In fact they obtain more precise results on functionnal inequalities like Prekopa-Leindler and Borell-Brascamp-Lieb inequalities.

Recently, there have been a lot of works on geometry of metric measure spaces. Lott-Villani  and Sturm have given independently a synthetic treatment of metric spaces having Ricci curvature bounded below by $k$ (see \cite{Villani,Sturm1,Sturm2}). All these works began by the result of precompactness of Gromov: the class of Riemannian manifolds of dimension $n$ and Ricci curvature bounded below by some constant $ k$ is precompact for a Gromov-Hausdorff metric. So the notion they develop for metric spaces has to generalize the one for Riemannian manifolds and has to be stable  by Gromov-Hausdorff convergence.  Their definition is about convexity properties of relative entropy on the Wasserstein space of probability  and is linked with optimal transportation.
Sturm in this context defines a Brunn Minkowski inequality with curvature $k$ (see \cite{Sturm2}).
 
The meaning of this inequality  may be not totally satisfactory. Indeed the inequality is depending on parameter $\Theta$ which equals $\inf_{k_\in K, l\in L} d(k,l) $ or $\sup_{k_\in K, l\in L} d(k,l) $ whether the curvature is positive (or null) or negative. It corresponds to the minimal or maximal length of geodesics between the two compacts $K$ and $L$.
%So the inequality can't be generalise to all spaces  
   However this is a direct implication from its dimension-curvature condition $CD(k,N)$ and this is this inequality that gives all the geometric consequences of their theory like for example a Bishop-Gromov theorem on the growth of balls.

There is another weak concept of curvature which is known as metric contraction property (see \cite {Ohta,Sturm2,Juillet}) and which  is implied by this Brunn-Minkowski inequality  at least in the case of curvature 0 and the $m\otimes m$ a.s. uniqueness of geodesics between two points of $X$.

As far as I know stability of Brunn-Minkowski inequality was not proven yet. This is the most interesting result we have in the paper (corollary   \ref{bmgeo}). For simplicity we will work only with the classical Brunn-Minkowski (i.e. with curvature 0) and explains how to extend our results in  the general case, with curvature $k$, in a remark. For doing this we introduce an approximated Brunn minkowski inequality since we need it during the proof. This fact is interesting in itself since it allows us to deal with discrete spaces.

In the second part of the paper we show that every metric measure space satisfying classical Brunn-Minkowski inequality can be approximated by discrete spaces with some approximated Brunn-Minkowski inequalities.

To avoid some problems between sets with zero measure we will work only with metric spaces $(X,d,m)$ where $(X,d)$ is Polish and $m$ a Borel measure on $(X,d)$ with full support, i.e. that charges every ball of $X$.

\section{Stability of Brunn-Minkowski inequality}
 
\begin{defi} \label{bm}
Given $h\geq 0$ and $N \in \N, N\geq 1$, we say that a  metric measure 
space $(X,d,\mu)$ satisfies the $h$ Brunn-Minkowski inequality of 
dimension $N$ denoted by $BM(N,h)$ if $\forall C_0,C_1 \subset X $ compacts, 
$\forall s\in [0,1]$, we have:
\begin{equation}
\mu^{1/N}({C_s}^h)\geq (1-s)\, \mu^{1/N}(C_0)  + s \,  \mu^{1/N}(C_1)
\end{equation} 
where
\begin{equation}
C_s^h=\left\{ x\in X / \exists (x_0,x_1)\in C_0\times C_1 /
\begin{array}{c c c}
|d(x_0,x)-s d(x_0,x_1)|&\leq & h \\ 
|d(x,x_1)-(1-s)d(x_0,x_1)|&\leq & h\\ 
\end{array}\right\}  
\end{equation}
\end{defi}

We call the set $C_s^h$ the set of $h$(-approximated) $s$-intermediate points between $C_0$ and $C_1$.
One can note that if $X$ is a geodesic space and $h=0$, it gives back 
the classical Brunn-Minkowski inequality for geodesic spaces. We shall often note $BM(N)$ instead of $BM(N,0)$. Another 
remark to be done is that this definition can be used for discrete spaces.\\
One can also note that if $X$ satisfy $BM(N,h)$ it will also satisfy $BM(N,h')$ for all $h'\geq h$.\\
In these notes we use the  following distance $\D$ between abstract metric measure spaces. We refer to \cite{Sturm1} for its properties.

\begin{defi}
 Let $(M,d,m)$ and $(M',d',m')$ be two metric measure spaces, their distance $\D$ is given by
$$\D((M,d,m),(M',d',m'))= \inf _{\hd,q} \left(\int_{M\times M'}\hd^2(x,x')dq(x,y)\right)^{1/2}
$$
where $\hd $ is a pseudo metric on $M\sqcup M'$ which coincides with $d$ on $M$ and with $d'$ on $M'$ and $q$ a coupling of the measures $m$ and $m'$. 
\end{defi}

\begin{theo} \label{bmth}
Let $(X_n,d_n,m_n)$ be a sequence of compact metric measure spaces which 
converges with respect to the distance $\D$ to another compact metric 
measure space $(X,d,m)$. If $(X_n,d_n,m_n)$ satisfies $BM(N,h_n)$ with 
$h_n \rightarrow h$ when $n$ goes to infinity, then $(X,d,m)$ satisfies 
$BM(N,h)$.  
\end{theo}
 In particular for compact geodesic spaces it implies directly the stability of 
the classical Brunn-Minkowski inequality with respect to the $\D$-convergence:\\ 

\begin{cor} \label{bmgeo}
Let $(X_n,d_n,m_n)$ be a sequence of compact geodesic spaces which 
converges with respect to the distance $\D$ to another compact metric 
measure space $(X,d,m)$, then $X$ is also a geodesic space. If $(X_n,d_n,m_n)$ satisfies $BM(N)$ then $(X,d,m)$ satisfies 
 also $BM(N)$.  
\end{cor}

We will make the proof of theorem \ref{bmth} only for  compact sets of strictly positive measure. 
The remarks after the proof will give the inequality for all mesurable sets.\\

 The idea of the proof is quite simple. We choose two compacts of the limit set $X$. Then we choose a good coupling of $X_n$ and $X$ and we construct two compacts of $X_n$ by dilating these compacts with respect to the pseudo-distance of the coupling and taking the restriction of this  two sets with $X_n$. The fact which makes things work is that the operation we did  doesn't lose two much measure. So, we can define a $s$-intermediate set in $X_n$ and apply  Brunn-Minkowski inequality in $X_n$. By the same construction as before, we construct a set in the limit set $X$ from the $s$-intermediate set in $X_n$ without loosing a lot of measure. To conclude we have to study the link between this set and  set of approximate $s$-intermediate points between initial compacts.\\

\textbf{Proof of Theorem \ref{bmth}}
Let $C_0,C_1$ two compacts of $X$ of strictly positive measure. Let $s\in [0,1]$. Choose $n$ so that 
$\D(X_n,X)\leq\frac{1}{2n}$. By definition of $\D$, there exists 
$\hat{d}$ a pseudo-metric on $X_n \sqcup X$ and $q$ a coupling of $m_n$ and 
$m$ so that 
$$\left( \int_{X_n\times X}\hd \acc 2 (x,y)dq(x,y)\right)\acc{1/2} \leq 
 \delta_n = \frac{1}{n}
$$

For $\en>0$ define $C_{n,i}\acc{\en} = \{ x \in X_n / \hd(x,C_i)\leq 
\en \}$ for $i=1,2$, these are compacts of $X_n$. They are indeed not empty for $n$ large enough and $\en$ well chosen, since being of strictly positive measure as we will see it.
We have 
\begin{eqnarray*}
m(C_0)&=& q(X_n \times C_0)\\
      &=& q(C_{n,0}\acc {\en} \times C_0) + q(\{ X_n\setminus 
C_{n,0} \acc {\en} \} \times C_0)
\end{eqnarray*}
But if $(x,y)\in \{X_n\setminus C_{n,0}^{\en} \} \times C_0$, then $\hd (x,y)\geq \en$, so 
\begin{eqnarray*}
q(\{ X_n\setminus C_{n,0} \acc {\en} \} \times C_0)&\leq & \int_{ \{X_n\setminus C_{n,0} \acc {\en} \} \times C_0 }
                                                                  \frac{\hd ^2(x,y)}{{\en}^2}dq(x,y)\\
                                                                  &\leq &  \frac{\delta_n^2}{{\en}^2} 
\end{eqnarray*}
which equals $\frac{1}{n}$ for $\delta_n=\frac{1}{n}$ and $\en=\frac{1}{\sqrt n }$.\\

On the other hand, we have:
\begin{eqnarray*}
m_n(C_{n,0}^{\en})&=& q(C_{n,0}^{\en} \times X)\\
                        &\geq& q(C_{n,0}\acc {\en} \times C_0)
\end{eqnarray*}
Consequently,
\begin{equation}\label{1}
m_n(C_{n,0}^\frac{1}{\sqrt n})\geq m(C_0)-\frac{1}{n}
\end{equation}
and identically
\begin{equation}\label{2}
m_n(C_{n,1}^\frac{1}{\sqrt n})\geq m(C_1)-\frac{1}{n}.
\end{equation}

Now consider the set $C_{n,s}^{\en,h_n}\subset X_n$ defined as in the definition (\ref{bm}) by 
$$C_{n,s}^{\en,h_n}= \left\{ x\in X_n / \exists (x_{n,0},x_{n,1})\in C_{n,0}^{\en}\times C_{n,1}^{\en} /
\begin{array}{c c c}
|d(x_{n,0},x)-s d(x_{n,0},x_{n,1})|&\leq & h_n \\ 
|d(x,x_{n,1})-(1-s)d(x_{n,0},x_{n,1})|&\leq & h_n\\ 
\end{array}\right\}  
$$
This is the set of all the $h_n$ $s$-intermediate points between $C_{n,0}^{\en}$ and $C_{n,1}^{\en}$. Since $X_n$ satisfies $BM(N,h_n)$, 
\begin{equation}\label{bmeq-n}
m_{n}^\frac{1}{n} (C_{n,s}^{\en,h_n})\geq (1-s)m_{n}^{1/N}(C_{n,0}^{\en})+s \, m_{n}^{1/N}(C_{n,1}^{\en})
\end{equation} 
We can now define $C_s^{\en, h_n} \subset X$ by 
$$C_{s}^{\en, h_n}=\{y\in X , \exists x \in C_{n,s}^{\en, h_n} \hd (x,y)\leq \en \}
$$
Similary to (\ref{1})  we have 
\begin{equation}\label{3}
m(C_{s}^{\frac{1}{\sqrt n},h_n})\geq m_n(C_{n,s}^{\en})-\frac{1}{n}
\end{equation} 
Now since $(x-\frac{1}{n})_+^{1/N}\geq x^{1/N}-(\frac{1}{n})^{1/N}$ for all $x\geq 0$, combining the inequalities (\ref{1}), (\ref{2}), (\ref{3}) and (\ref{bmeq-n}) give us, for $\en=\frac{1}{\sqrt n}$,
\begin{eqnarray*}
m^{1/N}(C_s^{\en, h_n}) &\geq & m_n^{1/N}(C_{n,s}^{\en, h_n})-(\frac{1}{n})^{1/N}\\
                                &\geq & (1-s)\,m_{n}^{1/N}(C_{n,0}^{\en})+s\, m_{n}^{1/N}(C_{n,1}^{\en})                                                -(\frac{1}{n})^{1/N}\\
                                &\geq & (1-s)\,m^{1/N}(C_0) +s \, m^{1/N}(C_1)-2(\frac{1}{n})^{1/N}	
\end{eqnarray*}

$C_s^{\en, h_n}$ is included in the set $K_s^{h_n + 4\en}$ of all the $h_n+4\en$ $s$-intermediate points between $C_0$ and $C_1$. Indeed, let $y\in C_s^{\en, h_n}$, by definition of this set, there exists $x\in C_{n,s}^{\en, h_n}$ so that $\hd (x,y)\leq \en$. By definition of $C_{n,s}^{\en, h_n}$, it follows that there exists $(x_{n,0},x_{n,1})\in C_{n,0}^{\en}\times C_{n,1}^{\en}$ satisfying 
$$\begin{array}{l c c}
|d_n(x,x_{n,0})- s\, d_n(x_{n,0},x_{n,1})|&\leq &h_n\\
|d_n(x,x_{n,1}) - (1-s)\, d_n(x_{n,0},x_{n,1})|&\leq &h_n.
\end{array}   
$$
There exists, by definition of $C_{n,i}^{\en}$ for  $i=1,2$, $(y_0,y_1) \in C_0\times C_1$ with $\hd (x_{n,0},y_0)\leq \en$ and  $\hd (x_{n,1},y_1)\leq \en$.
It follows:
\begin{eqnarray*}
|\hd(y,y_0)- s\, \hd(y_0,y_1)|&\leq & |\hd(y,y_0)-  \hd(x,x_{n,0})|+|\hd(x,x_{n,0})- s\, \hd(x_{n,0},x_{n,1})|\\
                                 &  &    + s\, |\hd(y_0,y_1)-  \hd(x_{n,0},x_{n,1})|\\
                              &\leq & h_n+ 4 \en .
\end{eqnarray*} 
and 
$$|\hd(y,y_1)- (1-s)\, \hd(y_0,y_1)|\leq h_n+ 4 \en .
$$ 

The sequence $(h_n+\en)_n$ is converging to $h$. We can extract a monotone sequence from it which will still be denoted by $h_n+\en$. There are two cases. The first one is when the extracting subsequence is non-decreasing. Then we have  $K_s^{h_n + 4\en} \subset K_s^{h}$. So, for all $n$, 
$$m^{1/N}(K_s^{h})\geq m^{1/N}(K_s^{h_n + 4\en})\geq (1-s)\,m^{1/N}(C_0) +s \, m^{1/N}(C_1)-2(\frac{1}{n})^{1/N}.$$
Letting $n$ goes to infinity gives the conclusion.\\
The second one, more interesting, is when the extracted subsequence is non-increasing. Then we have 
$$K_s^h= \bigcap_n K_s^{h_n + 4\en}.$$
Indeed if $y\in \bigcap_n K_s^{h_n + 4\en}$, for all $n\in \N$, $\exists (y_{n,0},y_{n,1})\in C_0\times C_1$ so that 
$$\begin{array}{l c c}
|d(y,y_{n,0})- s\, d(y_{n,0},y_{n,1})|&\leq &h_n+ 4\en\\
|d(y,y_{n,1})- (1-s)\, d(y_{n,0},y_{n,1})|&\leq &h_n+ 4\en .
\end{array}   
$$
By compactness of $C_0$ and $C_1$ we can extract another subsequence so that $y_{n,0}\rightarrow y_0 \in C_0$ and $y_{n,1}\rightarrow y_1 \in C_1$ and we have  
$$\begin{array}{l c c}
|d(y,y_{0})- s\, d(y_{0},y_{1})|&\leq &h\\
|d(y,y_{1})- (1-s)\, d(y_{0},y_{1})|&\leq &h
\end{array}   
.$$ The other inclusion is immediate.
This intersection is non-increasing so $$m^{1/N}(K_s^{h})= \lim_{n\rightarrow \infty} m^{1/N}(K_s^{h_n + 4\en})$$ which gives the conclusion
$$m^{1/N}(K_s^{h})\geq (1-s)\,m^{1/N}(C_0) +s \, m^{1/N}(C_1).$$

\textbf{Remark}
\begin{enumerate}

\item $BM(N)$ is directly implied by the condition $CD(O,N)$ of Sturm or Lott and Villani for the compact sets with a strictly positive measure (in fact for mesurable sets with strictly positive measure) (see \cite{Sturm2}). But if the measure $m$  is charging all the balls of the space and (if the space is geodesic), then the fact of having $BM(N)$ for all the compacts subspace with strictly positive measure implies $BM(N)$ for all compact subspaces.
Indeed if $(X,d,m)$ satifies $BM(N)$  for all the compact sets with a strictly positive measure and  if the measure $m$  is charging all the balls, if $C_0, C_1$ are compacts with $m(C_0)=0$ and $m(C_1)>0$ (the case $m(C_0)=m(C_1)=0$ is trivial) and $s\in[0,1]$. Define $C_0^\e=\{y\in X , \exists x\in C_0 / d(x,y)\leq \e\}$, $m(C_0^\e)>0$. Define $H_s^\e$ the set of all the $s$-intermediate points between $C_0^\e$ and $C_1$, By Brunn-Minkowski inequality we have:
$$m^{1/N}(H_s^\e)\geq (1-s)\, m^{1/N}(C_0^\e) + s\,m^{1/N}(C_1)\geq s\,m^{1/N}(C_1)$$
$H_s^\e$ is included in $K_s^{2\e}$ the set of all $2\e$ $s$-intermediate points between $C_0$ and $C_1$. 
As before $\bigcap_{\e>0} K_s^{2\e}$ is an non-increasing intersection equal to $K_s^0$ the set of all the exact $s$-intermediate points between $C_0$ and $C_1$. So 
$$m(K_s^0)=\lim_{\e \rightarrow 0}  K_s^{2\e}$$
which  gives the annonced result.
Consequently, on a metric measure space where the measure charges all the balls, $CD(0,N)$ implies $BM(N)$ for all compacts which in turns implies $MCP(0,N)$

\item In Polish spaces, Borel measures are regular which permits to pass from compact sets to measurable ones.  More precisely, if a Polish space satisfy $BM(N,h)$ for all his compact subsets, it also satisfies it for all his measurable subsets. Therefore, if the spaces $X_n$ and $X$ are only Polish (no more compacts),  the sets $C_{n,i}^{\en}$ for $i=1,2$ defined as above  may be no more compacts. However they will still be measurable since  closed, so  (\ref{bmeq-n})  will still stay true in this more general context. We can, consequently, drop the assumption of compactness of $X_n$ and $X$ in the theorem (\ref{bmth}) and its corollarry (\ref{bmgeo}).

\item We can do the same for the Brunn-Minkowski inequality with curvature $k$ by using  the definition given in \cite{Sturm2}. The only additional  thing to do is to control the parameter $\Theta$.
  But, with preceeding notations, we have $|\Theta(C_0,C_1) - \Theta(C_{n,0}^{\en},C_{n,1}^{\en})|\leq 2 \en $.

 \item We can prove also the same theorem for the multiplicative Brunn-Minkowski inequality (\ref{bminf}).
\end{enumerate}

\section{Discretizations of metric spaces}
Let $(M,d,m)$ be a given Polish measure space. For $h>0$, let $M_h=\{ x_i, i\geq 1\}$ be  a countable subspace of $M$ with $M=\bigcup_{i \geq 1} B_h(x_i)$. Choose $A_i\ \subset B_h(x_i), x_i \in A_i$ mutually disjoint and mesurable so that $\bigcup_{i\geq 1}A_i=M$. Consider the measure $m_h$ on $M_h$ given by $m_h(\{ x_i \})=m(A_i)$ for $i\geq 1$. We call $(M_h,d,m_h)$ a discretization of $(M,d,m)$.\\
It is proved in \cite{SturmB} that if $m(M)<\infty$ then 

$$(M_h,d,m_h)\xrightarrow[]{\D} (M,d,m).
$$

\begin{theo}\label{dth}
If $(M,d,m)$ satisfies $BM(N)$ then $(M_h,d,m_h)$ satisfies $BM(N,4h)$.
\end{theo}

 The proof is based on the two following facts.
 
 \begin{lemme}\label{dlem}
 
 \begin{enumerate}
 
 \item If $H \subset M_h$ then 
 \begin{equation}
 m(H^h) \geq m_h(H)
 \end{equation}
  where $H^h =\{x\in M , d(x,H)\leq h\}.$ 
  
 \item If $A\subset M$  mesurable and $A^h=\{x_i \in M_h, d(x_i,A)\leq h\}$ then
  \begin{equation}
  m_h(A^h)\geq m(A).
  \end{equation}  
 \end{enumerate}
 \end{lemme}
 
 \textbf{Proof of lemma \ref{dlem}}
 
 First, let $H\subset M_h$, we have
\begin{eqnarray*}
m_h(H)&=& \sum_{i/ x_i\in H}m(A_i)\\
       &=& m(\sqcup_{i / x_i \in H}A_i)\\
       &\leq & m(H^h)
\end{eqnarray*}
since $\sqcup_{i / x_i \in H}A_i \subset H^h= \{x\in M, \, d(x,H) \leq h \}$.\\
 For the second point, let $A \subset M$ mesurable, define $A^h$ as above, then
 
 \begin{eqnarray*}
m_h(A^h)&=& \sum_{i/ x_i\in A^h}m(A_i)\\
       &=& m(\sqcup_{i / x_i \in A^h}A_i)\\
       &\geq & m(A)
\end{eqnarray*}
  since $\sqcup_{i / x_i \in A^h}A_i \supset A$. Indeed if for some $j$, $A_j\cap A \neq \emptyset$ then there exists $a\in A$ with $d(x_j,a) \leq h$ so $x_j \in A^h$.\\

  \textbf{Proof ot theorem \ref{dth}}
  
  Let $H_0, H_1$ be two compacts of $M_h$ and $s\in [0,1]$. $H_0$ and $H_1$ consist of a finite or countable number of points $x_j$. Define $H_0 ^h,H_1^h \subset M$ by $H_i^h=\{x\in M, \exists x_j \in H_i / d(x_j,x)\leq h  \}$ for $i=1,2$. 
  By the first point of the lemma, for $i=1,2$
  \begin{equation}\label{1d}
  m(H_i^h) \geq m_h(H_i).
  \end{equation}
  
  Let $(H^h)_s \subset M$ be the set of all the $s$-intermediate points between $H_0^h$ and $H_1^h$ in the entire space $M$, i.e. 
  
  $$(H^h)_s=\left \{ x\in M , \exists (x_0,x_1) \in H_0^h \times H_1^h / \left|\begin{array}{l c c}
d(x,x_{0})&=& s\, d(x_{0},x_{1})\\
d(x,x_{1})&=& (1-s)\, d(x_{0},x_{1})
\end{array}\right. \right \}
  $$
  $BM(N)$ inequality on $M$ gives us 
  \begin{equation}\label{bmeqd}
m^{1/N} ((H^h)_s )\geq (1-s)m^{1/N}(H^h_0)+s\, m^{1/N}(H^h_1).
\end{equation} 
 
 As before by triangular inequality, we can see  $(H^h)_s$ is include in the set $\tilde{C}^{3h}_s$ of $3h$ $s$-intermediaire points in the whole space $M$ between $H_0$ and $H_1$. So the set $\tilde{H}^{4h}_s \subset M_h$  of $4h$ $s$-intermediate points between $H_0$ and $H_1$ in the discrete space $M_h$ contains the restriction at $M_h$ of the $h$ dilated of $(H^h)_s$. By the second point of the lemma we have 
 \begin{equation}\label{2d}
  m_h(\tilde{H}^{4h}_s) \geq m((H^h)_s).
 \end{equation} 
  
 Combining inequalities (\ref{1d}), (\ref{bmeqd}) and (\ref{2d}) ends the proof of the theorem.\\

  \textbf{Remark}
 If $(M,d,m)$ satisfies $BM(N,k)$ then $(M_h,d,m_h)$ satisfies $BM(N,k+4h)$.

  {\footnotesize %

}


\begin{thebibliography}{00}
\bibitem{Barthe} F. \textsc{Barthe}, \emph{Autour de l'in\'egalit\'e de {B}runn-{M}inkowski}. Ann. Fac. Sci. Toulouse Math. (6), (\textbf{2003}) vol 12, 27--178

\bibitem{SturmB}A.I. \textsc{Bonciocat} and K.T. \textsc{Sturm}, \emph{Mass transportation and rough curvature bounds for discrete spaces}. Preprint
  
 \bibitem{Burago} D. \textsc{Burago}, Y. \textsc{Burago} and S. \textsc{Ivanov}, \emph{A course in metric geometry}. Graduate Studies in Mathematics 33. American Mathematical Society, Providence, RI.(\textbf{2001})
  
 \bibitem{CorderoMCS} 
  D.\textsc{Cordero-Erausquin} , R. \textsc{McCann}  and
              M. \textsc{Schmuckenschl{\"a}ger} ,
     \emph{A {R}iemannian interpolation inequality \`a la {B}orell,
              {B}rascamp and {L}ieb}.
  Invent. Math. (\textbf{2001})
  vol 146, 219--257,

\bibitem{Gardner} R.J. \textsc{Gardner}, \emph{The Brunn-Minkowski inequality}. Bulletin of the American Mathematical Society  (\textbf{2001}) vol 39, n° 3, 355-405   

\bibitem{Juillet} N. Juillet \emph{Geometric Inequalities and Generalised Ricci Bounds in Heisenberg Group}, preprint 
    
\bibitem{Villani} J. \textsc{Lott} and C. \textsc{Villani}, \emph{Ricci curvature for metric-measure spaces via optimal transport}. Ann. of Math. (to appear).     
    
\bibitem{Ohta} S.I. \textsc{Ohta} \emph{On the measure contraction property of metric measure spaces} Comment. Math. Helv. (\textbf{2007}) Comment. Math. Helv. vol 82, 805--828 
\bibitem{Sturm1} K.T. \textsc{Sturm}, \emph{On the geometry of metric measure spaces. I}. Acta Math., in press.

\bibitem{Sturm2} K.T. \textsc{Sturm}, \emph{On the geometry of metric measure spaces. II}. Acta Math., in press. 

  
 
\bibitem{VillaniB} C. \textsc{Villani} \emph{Topics in optimal transportation}. Graduate Studies in Mathematics 58. American Mathematical Society (\textbf{2003}) 

\end{thebibliography}
\end{document}